\documentclass[12pt,reqno]{amsart}

\newcommand\version{May 16, 2010}

\usepackage{amsmath,amsfonts,amsthm,amssymb,amsxtra}

\newtheorem{theorem}{Theorem}

\newtheorem{lemma}[theorem]{Lemma}

\theoremstyle{definition}

\theoremstyle{remark}

\newcommand{\C}{\mathbb{C}}

\renewcommand{\epsilon}{\varepsilon}

\renewcommand{\phi}{\varphi}
\newcommand{\R}{\mathbb{R}}

\DeclareMathOperator{\im}{Im}

\DeclareMathOperator{\re}{Re}

\DeclareMathOperator{\sgn}{sgn}

\begin{document}

\title[Eigenvalue bounds --- \version]{Eigenvalue bounds for Schr\"odinger operators\\ with complex potentials}

\author{Rupert L. Frank}
\address{Rupert L. Frank, Department of Mathematics, Princeton University, Princeton, NJ 08544, USA}
\email{rlfrank@math.princeton.edu}

\thanks{\copyright\, 2010 by the author. This paper may be reproduced, in its entirety, for non-commercial purposes.}

\begin{abstract}
We show that the absolute values of non-positive eigenvalues of Schr\"o\-dinger operators with complex potentials can be bounded in terms of $L_p$-norms of the potential. This extends an inequality of Abramov, Aslanyan, and Davies to higher dimensions and proves a conjecture by Laptev and Safronov. Our main ingredient are the uniform Sobolev inequalities of Kenig, Ruiz, and Sogge.
\end{abstract}

\maketitle

\subsection*{Introduction and main result}

Recent years have seen an increasing interest in the spectral properties of non-selfadjoint differential operators. This theory is much less developed than the self-adjoint theory, and the absence of a spectral theorem and the lack of variational techniques poses great challenges, both on a theoretical and a numerical level. We refer to \cite{Da} for an overview, physical motivations, and a wealth of (counter)examples in the field.

In this paper we will be concerned with the situation of a non-selfadjoint Schr\"odinger operator $-\Delta+V$ in $\R^d$ with a complex-valued potential $V$ which decays at infinity (in some averaged sense). The essential spectrum of this operator acting in $L_2(\R^d)$ is $[0,\infty)$, and we are interested in the location of its discrete eigenvalues. It is a beautiful observation of \cite{AbAsDa} (see also \cite{DaNa}) that in dimension $d=1$, every eigenvalue $\lambda\in\C\setminus [0,\infty)$ of $-\frac{d^2}{dx^2}+V$ satisfies
\begin{equation}
 \label{eq:davies}
|\lambda|^{1/2} \leq \frac 12 \int_\R |V(x)| \,d x \,.
\end{equation}
In \cite{FrLaLiSe} we extended this bound to higher dimensions, to $L_p$-norms of $V$ with $p\neq 1$ and to sums of eigenvalues. The prize for these generalizations is, however, that the bounds are only valid for eigenvalues outside the cone $\{ z\in\C:\ |\im z| < \epsilon \re z \}$ around the positive half-axis, with constants blowing up as $\epsilon\to 0$. The recent papers \cite{LaSa,Sa} have led to a better understanding of eigenvalues lying close to the positive half-axis. The natural generalization of \eqref{eq:davies}, however,
\begin{equation}
 \label{eq:conj}
|\lambda|^\gamma \leq D_{\gamma,d} \int_{\R^d} |V(x)|^{\gamma+d/2} \,d x \,,
\end{equation}
for $d\geq 2$ and $0<\gamma\leq d/2$, was left as an open conjecture in \cite{LaSa}.

Our goal in this paper is to prove this conjecture for $0<\gamma\leq 1/2$. The motivation in \cite{LaSa} for the value $\gamma=d/2$ together with the example in \cite{IoJe} suggest that our range of $\gamma$'s is indeed best possible.

\begin{theorem}
 \label{main}
Let $d\geq 2$ and $0<\gamma\leq 1/2$. Then any eigenvalue $\lambda\in\C\setminus [0,\infty)$ of the Schr\"odinger operator $-\Delta +V$ satisfies
\begin{equation}
 \label{eq:main}
|\lambda|^\gamma \leq D_{\gamma,d} \int_{\R^d} |V(x)|^{\gamma+d/2} \,dx
\end{equation}
with a constant $D_{\gamma,d}$ independent of $V$.
\end{theorem}

We also prove the formal analogue of inequality \eqref{eq:main} for $\gamma=0$. For $d=3$ we can even evaluate the sharp value of the constant.

\begin{theorem}
 \label{main0}
Let $d\geq 3$. If
\begin{equation}
 \label{eq:main0}
D_{0,d} \int_{\R^d} |V(x)|^{d/2} \,dx < 1 \,,
\end{equation}
then the Schr\"odinger operator $-\Delta +V$ has no eigenvalue in $\C\setminus [0,\infty)$.
For $d=3$, one can take $D_{0,3}=\frac{4}{3^{3/2}\pi^{2}}$ and this is the smallest possible constant.
\end{theorem}

Actually, we shall prove that \eqref{eq:main} and \eqref{eq:main0} remain true if  $\int |V|^{\gamma+d/2} \,dx$ is replaced by the Morrey-Campanato norm
\begin{equation*}
\sup_{x,r} r^{d} \left( r^{-d} \int_{B_r(x)} |V(y)|^p \,dy \right)^{(2\gamma+d)/2p}
\end{equation*}
for $p_\gamma<p\leq \gamma+d/2$. This norm is weaker than the $L_{\gamma+d/2}$-norm entering in \eqref{eq:main} and \eqref{eq:main0} and even weaker than the weak-$L_{\gamma+d/2}$-norm. This allows for stronger local singularities than the standard $L_p$-norms (e.g., for $|x|^{-2}$ singularities in the estimate analogous to \eqref{eq:main0}). The resulting eigenvalue bounds seem to be new even for real-valued potentials.

\medskip

As is well-known, inequality \eqref{eq:main} is true for \emph{real-valued} potentials, even for all $\gamma>0$. This is an easy consequence of the Sobolev inequalities and the variational characterization of eigenvalues. There is no variational principle in the case of \emph{complex-valued} potentials, but this impass can be circumvented using the Birman-Schwinger principle similarly as in \cite{AbAsDa}. In contrast to the self-adjoint case, however, the classical Sobolev inequalities do not suffice for the proof of \eqref{eq:main}. We shall use the much deeper, `uniform Sobolev inequalities' of Kenig, Ruiz, and Sogge \cite{KeRuSo}.

A fundamental insight of Safronov \cite{Sa} is that methods from stationary scattering theory can be used to control eigenvalues of non-selfadjoint Schr\"odinger operators close to the essential spectrum. While the classical (self-adjoint) scattering theory of Agmon, Kato and Kuroda uses the Sobolev trace theorem for the restriction to the sphere of constant energy, the papers \cite{GoSc,IoSc} advocate the use of the Stein--Tomas restriction theorem and the closely related uniform Sobolev inequalities. In the same spirit, we will replace Sobolev's restriction theorem in \cite{Sa} by the uniform Sobolev inequalities in order to prove conjecture \eqref{eq:conj}.


\subsection*{Proofs}

The proofs of both theorems are based on the Birman--Schwinger principle: if $\lambda\in\C\setminus [0,\infty)$ is an eigenvalue of the Schr\"odinger operator $-\Delta +V$, then $-1$ is an eigenvalue of the Birman--Schwinger operator $V^{1/2} (-\Delta-\lambda)^{-1} |V|^{1/2}$. Here we use the abbreviation $V^{1/2} = (\sgn V) |V|^{1/2}$. If $-1$ is an eigenvalue of an operator, then the norm of this operator is at least $1$. Hence Theorem \ref{main} will follow if we can prove the bound
\begin{equation}
 \label{eq:equiv}
\left\| V^{1/2} (-\Delta-\lambda)^{-1} |V|^{1/2} \right\|^{\gamma+d/2} 
\leq D_{\gamma,d} |\lambda|^{-\gamma}  \int_{\R^d} |V(x)|^{\gamma+d/2} \,dx \,.
\end{equation}
Similarly, Theorem \ref{main0} will follow if we can prove the bound
\begin{equation}
 \label{eq:equiv0}
\left\| V^{1/2} (-\Delta+\lambda)^{-1} |V|^{1/2} \right\|^{d/2}
\leq D_{0,d}  \int_{\R^d} |V(x)|^{d/2} \,dx \,.
\end{equation}

As a warm-up, we begin with bound \eqref{eq:equiv0} for $d=3$. In this case, the integral kernel of $(-\Delta-\lambda)^{-1}$ is given by $(4\pi |x-y|)^{-1} e^{i\sqrt\lambda|x-y|}$, where we use the branch of the square root on $\C\setminus [0,\infty)$ with positive imaginary part. Hence this kernel is bounded in absolute value by $(4\pi |x-y|)^{-1}$, the Green's function. We conclude that for any $\phi,\psi\in L_2(\R^3)$,
\begin{align*}
 \left| \left(\phi, V^{1/2} (-\Delta-\lambda)^{-1} |V|^{1/2} \psi \right) \right| 
& \leq \frac1{4\pi} \iint \frac{|\phi(x)| |V(x)|^{1/2} |V(y)|^{1/2} |\psi(y)|}{|x-y|} \,dx\,dy \\
& \leq  \frac{2^{4/3}}{3\pi^{4/3}} \left\| \phi \,|V|^{1/2} \right\|_{6/5} \left\| \psi \,|V|^{1/2} \right\|_{6/5} \,,
\end{align*}
where in the last step we applied Lieb's sharp version of the Hardy--Littlewood--Sobolev inequality \cite{Li} (see also \cite{LiLo}). Finally, we bound using H\"older's inequality
$$
\left\| \psi \,|V|^{1/2} \right\|_{6/5} \leq \| V \|_{3/2}^{1/2} \, \|\psi\|_2 \,,
$$
and similarly for $\left\| \phi \,|V|^{1/2} \right\|_{6/5}$. Taking the supremum over all $\psi$ and $\phi$ with $L_2$-norm equal to one, we conclude that
$$
\left\| V^{1/2} (-\Delta+\lambda)^{-1} |V|^{1/2} \right\| \leq \frac{2^{4/3}}{3\pi^{4/3}} \| V \|_{3/2} \,.
$$
This proves the assertion of Theorem \ref{main0} for $d=3$ with the claimed constant.

This constant is sharp, even if the inequality is restricted to real-valued potentials. Indeed, the potential $V(x)=-3(1+|x|^2)^{-2}$ leads to a zero-energy resonance (with resonance function $\psi(x)=(1+|x|^2)^{-1/2}$), and any negative perturbation turns this resonance into a negative eigenvalue. Moreover, for this potential $\int |V|^{3/2}\,dx = 3^{3/2}\pi^2/4$.

We now turn to the general case. As we already mentioned, the key ingredient is the uniform Sobolev inequality by Kenig, Ruiz and Sogge \cite{KeRuSo}, which asserts that
\begin{equation}
 \label{eq:keruso}
\left\| (-\Delta-\lambda)^{-1} \right\|_{L_p\to L_{p'}} 
\leq C_{p,d} \, |\lambda|^{-\frac{d+2}{2}+ \frac dp}
\end{equation}
for $\frac{2d}{d+2} \leq p\leq \frac{2(d+1)}{d+3}$ if $d\geq 3$ and for $1< p\leq \frac{6}{5}$ if $d= 2$. Here, as usual, $\frac1p+\frac1{p'}=1$. To be more precise, the claimed $|\lambda|$-decay in \eqref{eq:keruso} is not stated in \cite{KeRuSo}, but once the inequality is established uniformly for $|\lambda|=1$, $\lambda\neq 1$, the decay follows easily by scaling. Moreover, in \cite{KeRuSo} the inequality is only proved for $d\geq 3$, but their argument extends to the two-dimensional case provided $p=1$ is excluded. (For $d=2$ one omits the gamma function factor in \cite[(2.17)]{KeRuSo} and verifies \cite[(2.20)]{KeRuSo} for $\re(\lambda)\in[-3/2,1)$.)

With \eqref{eq:keruso} at hand we can estimate similarly as before
\begin{align*}
 \left| \left(\phi, V^{1/2} (-\Delta+\lambda)^{-1} |V|^{1/2} \psi \right) \right| 
& \leq C_{p,d} \, |\lambda|^{-\frac{d+2}{2}+\frac dp} \left\| \phi \,|V|^{1/2} \right\|_{p} \left\| \psi \,|V|^{1/2} \right\|_{p} \\
& \leq C_{p,d} \, |\lambda|^{-\frac{d+2}{2}+\frac dp} \| V \|_{p/(2-p)} \|\phi\|_2 \|\psi\|_2 \,.
\end{align*}
 For given $0<\lambda\leq 1/2$ (as well as $\lambda=0$ for $d\geq 3$) we are allowed to choose $p=2(2\gamma+d)/(2\gamma+d+2)$. Taking again the supremum over all $\psi$ and $\phi$ with $L_2$-norm equal to one we arrive at \eqref{eq:equiv} and \eqref{eq:equiv0} with $D_{\gamma,d} = C_{p,d}^{\gamma+d/2}$. This completes the proofs of Theorems \ref{main} and \ref{main0}.


\subsection*{Extension to singular potentials}

We prove eigenvalue inequalities for potentials with stronger (local) singularities than those allowed in Theorems \ref{main} and \ref{main0}. For that purpose we measure the size of a potential by its norm in the Morrey-Campanato space $\mathcal L^{\alpha,p}(\R^d)$,
\begin{equation}
 \label{eq:morrey}
\|V\|_{\mathcal L^{\alpha,p}} := \sup_{x,r} r^{\alpha} \left( r^{-d} \int_{B_r(x)} |V(y)|^p \,dy \right)^{1/p} \,.
\end{equation}
One easily shows the (continuous) inclusions
$$
\mathcal L^{\alpha,p}(\R^d) \supset \mathcal L^{\alpha,r}(\R^d)  \supset L_{d/\alpha,\infty}(\R^d)  \supset \mathcal L^{\alpha,d/\alpha}(\R^d) = L_{d/\alpha}(\R^d)
\quad\text{for}\ 1\leq p<r<d/\alpha \,.
$$
(Indeed, one can prove that these inclusions are all strict.) We shall prove

\begin{theorem}
 \label{mainmorrey}
Let $d\geq 2$, $0<\gamma< 1/2$ and $(d-1)(2\gamma+d)/2(d-2\gamma)<p\leq \gamma+d/2$. Then any eigenvalue $\lambda\in\C\setminus [0,\infty)$ of the Schr\"odinger operator $-\Delta +V$ satisfies
\begin{equation}
 \label{eq:mainmorrey}
|\lambda|^\gamma \leq D_{\gamma,d,p} \, \sup_{x,r} r^{d} \left( r^{-d} \int_{B_r(x)} |V(y)|^p \,dy \right)^{(2\gamma+d)/2p} \,.
\end{equation}
Moreover, if $d\geq 3$, $(d-1)/2<p\leq d/2$ and
\begin{equation}
\label{eq:mainmorrey0}
D_{0,d,p} \, \sup_{x,r} r^{d} \left( r^{-d} \int_{B_r(x)} |V(y)|^p \,dy \right)^{d/2p} <1 \,,
\end{equation}
then the Schr\"odinger operator $-\Delta +V$ has no eigenvalue in $\C\setminus [0,\infty)$. 
\end{theorem}

For the proof we replace \eqref{eq:keruso} by a weighted $L_2$-estimate. The $L_2$-space on $\R^d$ with measure $\omega(x)dx$ is denoted by $L_2(\omega)$.

\begin{lemma}\label{chsa}
 Let $4/3< \alpha< 2$ if $d=2$, $2d/(d+1)<\alpha\leq 2$ if $d\geq 3$ and let $(d-1)/2(\alpha-1)<p\leq d/\alpha$. Then for all $0<\omega\in \mathcal L^{\alpha,p}(\R^d)$,
\begin{equation}
 \label{eq:chsa}
\left\| (-\Delta-\lambda)^{-1} \right\|_{L_2(\omega^{-1})\to L_{2}(\omega)} 
\leq C_{d,\alpha,p}   \|\omega\|_{\mathcal L^{\alpha,p}} \, |\lambda|^{-1+\alpha/2} \,.
\end{equation}
\end{lemma}

This bound is essentially contained in the works \cite{ChSa, ChRu}, but not in this explicit form and only for $d\geq 3$. Therefore we provide a short sketch of a proof in the final section below. Assuming it for the moment, we choose a strictly positive function $g\in \mathcal L^{\alpha,p}$ (a Gaussian, say) and define $V_\epsilon(x) := \sup\{|V(x)|,\epsilon g(x)\}$. Then by Lemma \ref{chsa}
\begin{align*}
 \left| \left(\phi, V^{1/2} (-\Delta+\lambda)^{-1} |V|^{1/2} \psi \right) \right| 
& \leq C_{d,\alpha,p} |\lambda|^{-1+\alpha/2} \| V_\epsilon \|_{\mathcal L^{\alpha,p}} \,
 \left\| \sqrt{|V|/ V_\epsilon}\, \phi \right\|_{2} 
 \left\| \sqrt{|V|/ V_\epsilon}\, \psi \right\|_{2} \\
& \leq C_{d,\alpha,p} |\lambda|^{-1+\alpha/2} \| V_\epsilon \|_{\mathcal L^{\alpha,p}} \,
 \left\| \phi \right\|_{2} 
 \left\| \psi \right\|_{2} \,.
\end{align*}
Hence, as $\epsilon\to 0$,
$$
\left\| V^{1/2} (-\Delta+\lambda)^{-1} |V|^{1/2} \right\| 
\leq C_{d,\alpha,p} |\lambda|^{-1+\alpha/2} \| V \|_{\mathcal L^{\alpha,p}} \,,
$$
which is the analogue of \eqref{eq:equiv} and \eqref{eq:equiv0}. Theorem \ref{mainmorrey} follows by choosing $\alpha=2d/(2\gamma+d)$.


\subsection*{Extension to slowly decaying potentials}

In our previous theorems, we essentially required the potential to decay like $|x|^{-\rho}$ with $\rho>2d/(d+1)$. In order to deal with the case $\rho>1$, we recall the following weighted norm estimate for the resolvent $(-\Delta-\lambda)^{-1}$, which plays a crucial role in stationary scattering theory,
\begin{equation}
 \label{eq:agmon}
\left\| (-\Delta-\lambda)^{-1} \right\|_{L_2((1+|x|^2)^\alpha) \to L_2((1+|x|^2)^{-\alpha})}
\leq C_{d,\alpha} \, |\lambda|^{-1/2}
\qquad\text{if}\ \alpha>1/2 \,.
\end{equation}
The dependence on $\lambda$ can be improved to $(1+|\lambda|)^{-1/2}$, but \eqref{eq:agmon} is enough for our purposes. We refer to \cite{Sy} for the $\alpha$ dependence. In contrast to \eqref{eq:keruso} and \eqref{eq:chsa}, $|\lambda|$ enters \eqref{eq:agmon} with the (optimal) power $-1/2$. Improvements in this endpoint case can be found in \cite{BaRuVe}.

The same argument as in the previous subsection, but with $V_\epsilon$ replaced by $(1+|x|^2)^\alpha$, leads to the bound
$$
|\lambda|^{1/2} \leq C_{d,\alpha} \, \sup_x (1+|x|^2)^\alpha |V(x)| \,,
\qquad\alpha>1/2 \,,
$$
thus rederiving the result from \cite{Sa}. By interpolating this with Theorem \ref{main} we obtain 

\begin{theorem}
Let $d\geq 2$, $\gamma>1/2$ and $\alpha>\gamma-1/2$. Then any non-positive eigenvalue $\lambda\in\C\setminus [0,\infty)$ of the Schr\"odinger operator $-\Delta +V$ satisfies
\begin{equation}
 \label{eq:interpol} 
|\lambda|^\gamma \leq C_{d,\gamma,\alpha} \int_{\R^d} |V(x)|^{2\gamma+(d-1)/2} (1+|x|^2)^\alpha \,dx \,.
\end{equation}
\end{theorem}

Indeed, interpolation between \eqref{eq:keruso} with $p_d:=2(d+1)/(d+3)$ and \eqref{eq:agmon} yields
\begin{equation*}
\left\| (1+|x|^2)^{-\alpha\theta/2} (-\Delta-\lambda)^{-1} (1+|x|^2)^{-\alpha\theta/2} \right\|_{L_p \to L_{p'}}
\leq C_{d,\alpha,p} \, \, |\lambda|^{-\tfrac{1-\theta}{d+1}-\tfrac\theta2} \,,
\quad \tfrac1p=\tfrac{1-\theta}{p_d} + \tfrac\theta2 \,,
\end{equation*}
for $\alpha>1/2$, and then the argument of the previous subsection implies \eqref{eq:interpol}.


\subsection*{Proof of Lemma \ref{chsa}}
 As a preliminary step we deduce from the works \cite{ChWh,SaWh} that the operator $I_\alpha$ of convolution with $|x|^{-d+\alpha}$, $0<\alpha<d$, satisfies
\begin{equation}
 \label{eq:chwh}
\| I_\alpha \|_{L_2(w^{-1})\to L_2(w)} \leq C_{d,\alpha,s} \|w\|_{\mathcal L^{\alpha,s}}
\end{equation}
for any $1<s\leq d/\alpha$. Indeed, the constant $C_r$ in \cite[(1.10)]{SaWh} for $w=v^{-1}$ and $p=2$ is bounded by a constant times $\|w\|_{\mathcal L^{\alpha,r}}$, and this constant enters multiplicatively in \cite[(2.9)]{SaWh}.

After this preliminary step we turn to the proof of estimate \eqref{eq:chsa}. Following \cite{KeRuSo,ChSa} we consider the analytic family of operators $\omega^{\zeta/2}(-\Delta-\lambda)^{-\zeta}\omega^{\zeta/2}$ with $0\leq\re\zeta\leq (d+1)/2$. The bounds in \cite{KeRuSo} for the (explicit) integral kernel of $(-\Delta-\lambda)^{-\zeta}$ show that
$$
\left| (-\Delta-\lambda)^{-\zeta}(x,x') \right| \leq C_d e^{C|\im\zeta|} |x-x'|^{-(d+1)/2+\re\zeta}
$$
for $(d-1)/2\leq\re\zeta\leq (d+1)/2$ with $C$ independent of $\lambda$ satisfying $|\lambda|=1$. This, together with \eqref{eq:chwh}, proves that
\begin{align*}
 \| \omega^{\zeta/2}(-\Delta-\lambda)^{-\zeta}\omega^{\zeta/2} \|_{L_2\to L_2} 
& \leq C e^{C|\im\zeta|} \|\omega^{\re\zeta} \|_{\mathcal L^{\re\zeta+(d-1)/2,s}} \\
& = C e^{C|\im\zeta|} \|\omega \|_{\mathcal L^{1+(d-1)/(2\re\zeta),s\re\zeta}}^{\re\zeta}
\end{align*}
for $s>1$. Complex interpolation with the trivial bound for $\re\zeta=0$ yields
$$
\| \omega^{1/2}(-\Delta-\lambda)^{-1}\omega^{1/2} \|_{L_2\to L_2}
\leq \tilde C \|\omega \|_{\mathcal L^{1+(d-1)/(2\re\zeta),s\re\zeta}}
$$
for $(d-1)/2\leq\re\zeta\leq (d+1)/2$ if $d\geq 3$ and for $1\leq \re\zeta\leq 3/2$.
This is the claimed bound for $|\lambda|=1$, and the case of general $\lambda$ follows by scaling. (Note that $2d/(d+1)\leq\alpha=1+(d-1)/(2\re\zeta)\leq 2$ corresponds to $(d-1)/2\leq\re\zeta\leq (d+1)/2$.)

It remains to prove the case $\alpha\in(3/2,2)$ for $d=2$. A similar (but simpler) analysis as in \cite{KeRuSo} shows that for $d=2$ and any $3/2<\alpha<2$ one has
$$
\left| (-\Delta-\lambda)^{-1}(x,x') \right| \leq C_\alpha |x-x'|^{-2+\alpha}
$$
uniformly in $|\lambda|=1$. This, together with \eqref{eq:chwh}, completes the proof of Lemma \ref{chsa}.


\bigskip

The author wishes to thank A. Laptev and O. Safronov for useful correspondence.


\bibliographystyle{amsalpha}

\end{document}